\documentclass[11pt,leqno]{article}

\usepackage{amsmath,amsthm}
\usepackage {latexsym}
\usepackage{amssymb}

\newcommand{\beeq}{\begin{equation}}
\newcommand{\eneq}{\end{equation}}

\newcommand{\bear}{\begin{eqnarray}}
\newcommand{\eear}{\end{eqnarray}}

\newcommand{\supp}{\mbox{\rm supp}}

\newcommand{\half}{\frac{1}{2}}

\newcommand{\eps}{{\varepsilon}}
\newcommand{\R}{{\mathbb R}}

\newcommand{\Z}{{\mathbb Z}}

\newcommand{\Compl}{{\mathbb C}}
\newcommand{\Half}{{\mathbb H}}

\newcommand{\calS}{{\mathcal S}}

\newcommand{\les}{\lesssim}

\newcommand{\etatil}{\tilde{\eta}}
\newcommand{\rtil}{\tilde{r}}
\newcommand{\omtil}{\tilde{\omega}}

\newcommand{\Laplace}{\triangle}

\newcommand{\la}{\lambda}

\def\Lap{\Delta}
\def\nn{\nonumber}
\def\Sal{S_\alpha^{(R)}}

\def\gal{g_\alpha}
\def\fal{f_\alpha}

\newtheorem{theorem}{Theorem}
\newtheorem{lemma}[theorem]{Lemma}

\newtheorem{cor}[theorem]{Corollary}
\newtheorem{prop}[theorem]{Proposition}
\newtheorem{proposition}[theorem]{Proposition}

\theoremstyle{remark}

\def\la{\langle}
\def\ra{\rangle}

\def\norm[#1][#2]{\Vert #1 \Vert_{#2}}

\def\xial{\xi_\alpha}
\def\rank{{\rm rank}}
\def\slashint{{\int\!\!\!\!\!\!-\!\!\!}}

\begin{document}

\title{A Limiting Absorption Principle for the three-dimensional Schr\"odinger
equation with $L^p$ potentials}

\author{M.\ Goldberg, W.\ Schlag}
\date{}

\maketitle

\section{Introduction}

Agmon's fundamental work~\cite{Ag} establishes the bound, known as the limiting absorption principle,
\beeq
\label{eq:LA}
\sup_{\lambda>\lambda_0,\,\eps>0}\big\|\big(-\Laplace+V-(\lambda^2+i\eps)\big)^{-1}\big\|_{L^{2,\sigma}(\R^d)\to L^{2,-\sigma}(\R^d)}
<\infty
\eneq
provided that $\lambda_0>0$, $(1+|x|)^{1+}|V(x)|\in L^\infty$ and $\sigma>\half$. Here
\[L^{2,\sigma}(\R^d)=\{(1+|x|)^{-\sigma}\,f\::\: f\in L^2(\R^d)\}\]
is the usual weighted $L^2$. The bound~\eqref{eq:LA} is obtained from the same estimate for $V=0$
by means of the resolvent identity. This bound for the free resolvent is related to the
so called trace lemma, which refers to the statement that for every $f\in L^{2,\half+}$
there is a restriction of $\hat{f}$ to any (compact) hypersurface, and this restriction belongs
to $L^2$ relative to surface measure. Note that this fact does not require any curvature
properties of the hypersurface - in fact, it is proved by reduction to flat surfaces.
Another fundamental restriction theorem is the Stein-Tomas theorem, see~\cite{stein}.
It requires the hypersurfaces $\calS\subset \R^d$ with $d\ge2$ to
have non vanishing Gaussian curvature, and states that
\beeq
\label{eq:st_tom}
\int_{\calS} |\hat{f}(\omega)|^2\, \sigma(d\omega) \le C\|f\|^2_{L^p(\R^d)} \text{\ \ where\ \ }p=\frac{2d+2}{d+3}.
\eneq
It is not hard to see that the related estimate for the free resolvent in $\R^3$ is given by
\beeq
\label{eq:free_res}
 \|R_0(\lambda^2+i0)\|_{\frac43\to 4} = C\, \lambda^{-\half}\text{\ \ for\ \ }\lambda>0.
\eneq
This fact depends on the oscillation in the resolvent, i.e., on the exponential in
\beeq
\label{eq:res_ker}
 R_0(\lambda^2+i0)(x,y) = \frac{e^{i\lambda|x-y|}}{4\pi |x-y|}.
\eneq
In contrast, using the denominator alone one obtains that \beeq
\label{eq:useless}
\sup_{\lambda}\|R_0(\lambda^2+i0)\|_{\frac65\to 6}\le C
\eneq
via fractional integration.
In analogy with Agmon's work, it is natural to ask for which potentials~\eqref{eq:free_res} can be extended
to the perturbed operators $H=-\Laplace+V$.
In this paper we show that this is the case
for real-valued $V\in L^p(\R^3)\cap L^{\frac32}(\R^3), p > \frac32$,
and suggest two possible extensions.

\begin{theorem}
\label{thm:agmon} Let $V\in L^p(\R^3)\cap L^{\frac32}(\R^3), p > \frac32$
be real-valued.  Then for every $\lambda_0>0$, one has
\begin{equation}
\label{eq:V34}
\sup_{0<\eps<1,\;\lambda\ge\lambda_0}
\Big\|(-\Laplace+V - (\lambda^2+i\eps))^{-1}\Big\|_{\frac43\to 4} \le C(\lambda_0,V)\;\lambda^{-\half}.
\end{equation}
In particular, the spectrum of $-\Laplace+V$ is purely absolutely continuous on $(0,\infty)$.
\end{theorem}

We  also formulate dynamical consequences of this result, in
particular the existence and completeness of the wave operators.
This theorem is the analogue of the classical Kato-Agmon-Kuroda
theorem, see \cite{RS4}, Theorem~XIII.33. It of course requires
the absence of imbedded eigenvalues. In the classical context one
uses Kato's theorem for that purpose. Here we wish to use a result on the
absence of imbedded eigenvalues that only requires an
integrability condition on~$V$. One such result was obtained by
Ionescu and Jersion~\cite{IJ}, namely:

\begin{theorem} \label{thm:IJ}
Let $V\in L^{\frac32}(\R^3)$. Suppose $u\in W^{1,2}_{\rm loc}(\R^3)$ satisfies
$(-\Laplace +V)u=\lambda^2 u$ where $\lambda\ne0$ in the sense of
distributions. If, moreover,
$\|(1+|x|)^{\delta-\half}u\|_2<\infty$ for some $\delta>0$, then $u\equiv 0$.
\end{theorem}

The weighted $L^2$-condition with $\delta>0$ is natural in view of the Fourier
transform of the surface measure of $S^2$,
which is a generalized eigenfunction
of the free case and decays like $(1+|x|)^{-1}$.
As far as local regularity of the potential is concerned, the requirement that
$V \in L^{3/2}_{\rm loc}$ is essentially optimal.  There exist examples of
$V \in L^{3/2}_{\rm weak}$ for which $-\Lap + V$ admits compactly supported
eigenfunctions \cite{KocTat}.  The necessary decay condition
on $V$ is less clearly delineated:
Ionescu and Jerison found a smooth real-valued potential $V$ which
lies in $L^q(\R^3)$ for all $q>2$ but such that for $-\Laplace+V$ imbedded
eigenvalues exist. Their example decays like $r^{-1}$ in some directions,
and like $r^{-2}$ in other directions.
They further conjectured that their main result (Theorem 2.1 in \cite{IJ})
remains valid for potentials $V\in L^2(\R^3)$.
Recent work by Koch and Tataru appears to verify this conjecture
\cite{KocTat2}, and futher refinements which allow potentials to exhibit both
$L^{3/2}_{\rm loc}$ singularities and $L^2$ decay seem possible as well.
The proof of any such conjecture would immediately increase the scope of
Theorem~\ref{thm:agmon}, as described below.
\begin{prop} \label{prop:extension}
The following inferences are valid:

\noindent
1.  If the conclusion of Theorem~\ref{thm:IJ} holds for all $V \in L^p(\R^3),
\frac32 \le p < 2$, as is suggested by \cite{KocTat2}, then the conclusion of
Theorem~\ref{thm:agmon} also holds for all $V\in L^p(R^3)$.

\noindent
2. More generally,
if the conclusion of Theorem~\ref{thm:IJ} holds for some $V\in L^p(\R^3) +
L^q(\R^3), \frac32 < p,q < 2$, then the conclusion of Theorem~\ref{thm:agmon}
also holds for this $V$.
\end{prop}

By Kato's theory of $H$-smoothing operators, see~\cite{Kato}, it is well-known that
the limiting absorption principle for the resolvent gives rise to estimates for the
evolution $e^{itH}$ known as smoothing estimates.
This is a much studied class of bounds, see~\cite{SJ}, \cite{VE}, \cite{CS1}, \cite{CS2},
\cite{BAK}, \cite{Doi}, \cite{simon}. 
In fact, the Fourier transform establishes a
link between the resolvent and the evolution that in a  precise sense allows one to
state that a certain class of estimates  on the evolution is equivalent to corresponding ones for the
resolvent, see~\cite{Kato}. In the free case, the $\frac43\to 4$ bound for the resolvent 
corresponds to the following smoothing bound for the free evolution:
\[ \sup_{\|F\|_4\le 1} \int_{-\infty}^\infty \Big\|F(-\Laplace)^{\frac18}
e^{it\Laplace} f\Big\|_2^2\, dt \le C\,\|f\|_2^2.\]
However, this bound is known, see the work of Ruiz and Vega \cite{RuiVeg}.
For the perturbed evolution, $H = -\Lap + V$, one can prove
similar estimates by means of Theorem~\ref{thm:agmon}, but we do not pursue this here.
See the work of Ionescu and the second author~\cite{IonSch} for statements of this type.

This paper is organized as follows: In Section~\ref{sec:free}
we prove the bounds on the free resolvent that are needed in order to
prove Theorem~\ref{thm:agmon}. Our main new bounds involve $R_0(\lambda^2+i0)$
acting on functions whose Fourier transform vanish on $\lambda S^2$.
In Section~\ref{sec:perturb} we apply these bounds in the context of the usual
resolvent identity/Fredholm alternative type arguments to deal with
$-\Laplace+V$.  This of course requires Theorem~\ref{thm:IJ}.
Finally, in Section~\ref{sec:more} we return to the
free resolvent and prove some end point results.

\section{The free resolvent}\label{sec:free}

This section develops some estimates on the free resolvent given by~\eqref{eq:res_ker}.
These estimates are motivated on the one hand by the Stein-Tomas theorem~\eqref{eq:st_tom}, and
on the other hand, by the applications to the perturbed operator $H=-\Laplace+V$, see Theorem~\ref{thm:IJ}. For what follows, it will be helpful to keep in mind that for real $\lambda$,
\[ [R_0(\lambda^2+i0)-R_0(\lambda^2-i0)]f = C(\lambda)\cdot(\widehat{\sigma_{\lambda S^2}}\ast f),\]
which is exactly of the form $T^*T$, $T$ being the restriction operator to the sphere $\lambda S^2$.
Thus $T^* T: L^{\frac43}(\R^3)\to L^4(\R^3)$ in view of~\eqref{eq:st_tom}.

We will denote by $\Half$ the closed upper half-plane in $\Compl$,
and state most of
our results for $\lambda \in \Half$.  For any positive real number $\lambda$,
we have the boundary identites
$$(\lambda + i0)^2 = \lambda^2 + i0 \quad {\rm and} \quad
  (-\lambda + i0)^2 = \lambda^2 - i0,$$
therefore estimates which hold uniformly out to $\partial\Half$ are of
particular importance.

\begin{lemma}
\label{lem:free}
Let $\lambda \in \Half$ be any nonzero element, and $p=\frac43$. Then
$R_0(\lambda^2) : L^p(\R^3)\to L^{p'}(\R^3)$, with operator norm bounded by
$|\lambda|^{-\frac12}$.
\end{lemma}
As suggested above, the proof follows a complex-interpolation argument
strongly reminiscent of the proof of \eqref{eq:st_tom}.  For full details see
Theorem~2.3 in \cite{KRS}, which establishes this bound for a more general
family of inverses of second-order differential operators.

\begin{lemma}\label{lem:free2}
Let $\lambda \in \Half$ be any nonzero element.
For each pair of exponents $1 < p \le \frac43$, $3p \le q \le \frac{3p}{3-2p}$
there exist constants $C_{p,q} < \infty$
such that
$$\norm[R_0(\lambda^2)f][L^{q}] \le C_{p,q} |\lambda|^{3/p - 3/q- 2}
\norm[f][L^{p}]$$
For each exponent $\frac43 \le p < \frac32$, $\frac{p}{3-2p} \le q \le
\frac{3p}{3-2p}$ there exist constants
$C_{p,q} < \infty$ such that
$$\norm[R_0(\lambda^2)f][L^{p*}] \le C_{p,q}|\lambda|^{3/p - 3/q - 2}
\norm[f][L^p]$$
\end{lemma}
\begin{proof}
The case $p = \frac43, q = 4$ is Lemma~\ref{lem:free} above.
Since $R_0(\lambda^2)$ is realized as a convolution with a kernel
satisfying $|K_\lambda(x)| \le |4\pi x|^{-1}$, the cases $q = \frac{3p}{3-2p},
1 < p < \frac32$ are precisely the Hardy--Littlewood--Sobolev inequality.
Note that the scaling exponent for $\lambda$ is zero for these pairs $(p,q)$.
All intermediate cases $(p,q)$ then follow by interpolation.
At the endpoint $p=1, q=3$, we see that
$R_0^\pm(\lambda^2)$ maps $L^1(\R^3)$ to weak-$L^3(\R^3)$ uniformly in
$\lambda$, by considering the norm
$$\norm[f][L^3_{\rm weak}(\R^3)] = \sup_{A\subset\R^3,\, |A|< \infty}
|A|^{-\frac23}\int_A |f(x)|\,dx$$
which is equivalent to the usual weak-$L^3$ ``norm'' and satisfies a
triangle inequality, see Lieb, Loss~\cite{LL}, Section~4.3
The cases $q = 3p$, $1< p < \frac43$ follow
by Marcinkiewicz interpolation, and $q = \frac{p}{3-2p}$,
$\frac43 < p < \frac32$ by duality.
\end{proof}

The following results deal with functions whose Fourier transform vanishes on
$S^2$. The first lemma yields a H\"older bound for the $L^2$ norms of the
restrictions to spheres close to~$S^2$.

\begin{lemma}
Let $1\le p<\frac43$ and set $\gamma=\frac{2}{p}-\frac32$. Then for all $|\delta|<\half$ one has
\begin{equation}
\label{eq:hold}
\|\hat{f}((1+\delta)\cdot)\|_{L^2(S^2)} \les |\delta|^{\gamma} \|f\|_{L^p(\R^3)}
\end{equation}
for all $f\in L^p(\R^3)$ with $\hat{f}=0$ on $S^2$.
\end{lemma}
\begin{proof}
Let $\sigma_{(1+\delta)S^2}$ be the normalized measure on $(1+\delta)S^2$. Then one has
\bear
\|\hat{f}((1+\delta)\cdot)\|_{L^2(S^2)}^2 &=& \langle f\ast \widehat{\sigma_{(1+\delta)S^2}},f \rangle =
\langle f\ast [\widehat{\sigma_{(1+\delta)S^2}}-\widehat{\sigma_{S^2}}],f \rangle \nn \\
&=& \sum_{j=0}^\infty \langle f\ast K_j, f\rangle \nn
\eear
where $K_j(x) = \Big( \widehat{\sigma_{(1+\delta)S^2}}-\widehat{\sigma_{S^2}}\Big)\chi_j$
and $\{\chi_j\}_{j\ge0}$ are a standard dyadic partition of unity. Since
$\|\widehat{\sigma_{(1+\delta)S^2}}-\widehat{\sigma_{S^2}}\|_\infty\les \delta$, it follows that
\begin{equation}
\|K_j\|_\infty \les \left\{
\begin{array}{cc} \delta & \text{if\ \ }2^j<\delta^{-1} \\
                  2^{-j} & \text{if\ \ }2^j\ge \delta^{-1}
\end{array}
\right.\nn
\end{equation}
Thus $\|K_j\|_\infty \les \min(\delta,2^{-j}):=\alpha_j$.
Moreover,
\bear
\|\widehat{K_j}\|_\infty &=& \big\|\big(\sigma_{(1+\delta)S^2}-\sigma_{S^2}\big)\ast \widehat{\chi_j}\big\|_\infty
\nn \\
&=& \Big| \int \widehat{\chi_j}(\xi-\eta)\,\sigma_{(1+\delta)S^2}(d\eta)-\int \widehat{\chi_j}(\xi-\eta)\,\sigma_{S^2}(d\eta)\Big| \nn \\
&=& \Big| \int \Big[\widehat{\chi_j}(\xi-(1+\delta)\eta)-\widehat{\chi_j}(\xi-\eta)\Big]\,\sigma_{S^2}(d\eta)\Big| \nn \\
&\les& \min(2^{2j}\delta,2^j):=\beta_j. \nn
\eear
If $1 < p <\frac43$, let $\frac{1}{p}=\frac{\theta}{1}+\frac{1-\theta}{2}$
so that $\theta>\half$.
Then $\|K_j\ast f\|_{p'}\les \alpha_j^\theta\beta_j^{1-\theta}\|f\|_p$ for
all $j\ge0$. Summing over $j$ yields the desired bound.  In the case $p = 1$,
the estimate
$\norm[\widehat{\sigma_{(1+\delta)S^2}}-\widehat{\sigma_{S^2}}][\infty]
\les \delta$ mentioned above suffices to show that
$\norm[\hat{f}((1+\delta)\cdot)][L^2(S^2)] \les \delta^{\frac12}$.

\end{proof}

The point of the following proposition is that one can take $\delta>0$
in~\eqref{eq:main_est}. In the following section, this will allow us to apply
Theorem~\ref{thm:IJ}.

\begin{prop}
\label{prop:ess}
Let $1\le p<\frac43$. Then for any $\delta<\frac12-\frac{2}{p'}$ one has
\begin{equation}
\label{eq:main_est}
\sup_{\eps>0}\Big\|(1+|x|)^{\delta-\half} R_0(1\pm i\eps)f\Big\|_2 \les \|f\|_{p}
\end{equation}
for any $f\in L^p(\R^3)$ so that $\hat{f}=0$ on $S^2$.
\end{prop}
\begin{proof} We first consider the case where
\begin{equation}
\label{eq:supp_res}
\supp(\hat{f})\subset \{\xi\in\R^3\::\:\half<|\xi|<2\}.
\end{equation}
Let $\chi$ be a smooth, radial, bump function around zero so that
$\hat{\chi}$ is compactly supported. Let $R\gg1$.
Then
\bear
\|\chi(\frac{\cdot}{R})R_0(1+i\eps)f\|_2^2 &=& R^6 \int_{\R^3}\int_{\R^3}
\hat{\chi}(R(\xi-\eta))\frac{\hat{f}(\eta)}{|\eta|^2-1-i\eps}\,d\eta
\int_{\R^3}
\hat{\chi}(R(\xi-\etatil))\frac{\overline{\hat{f}(\etatil)}}{|\etatil|^2-1+i\eps}\,d\etatil \nn \\
&=& R^3 \int_{\R^6} \rho(R(\eta-\etatil))\, \frac{\hat{f}(\eta)}{|\eta|^2-1-i\eps}
\frac{\overline{\hat{f}(\etatil)}}{|\etatil|^2-1+i\eps} \, d\eta\,d\etatil, \label{eq:odin}
\eear
where we have set
\bear
\int_{\R^3}\hat{\chi}(R(\xi-\eta))\hat{\chi}(R(\xi-\etatil)) \,d\xi &=&
R^{-3}\int_{\R^3}\hat{\chi}(\zeta-R\eta)\hat{\chi}(\zeta-R\etatil) \,d\zeta \nn \\
&=& R^{-3}\int_{\R^3}\hat{\chi}(\zeta-R(\eta-\etatil))\hat{\chi}(\zeta) \,d\zeta \nn \\
&=:& R^{-3} \rho(R(\eta-\etatil)). \nn
\eear
Note that $\rho$ is a compactly supported smooth bump-function. Introducing
polar coordinates in~\eqref{eq:odin} yields uniformly in~$\eps\ne0$  (recall~\eqref{eq:supp_res})
\bear
\eqref{eq:odin} &=& R^3\int_{\R^3}\int_0^\infty\int_{S^2} \rho(R(\eta-\rtil\omtil)) \frac{\hat{f}(\eta)}{|\eta|^2-1-i\eps}
\frac{\overline{\hat{f}(\rtil\omtil)}}{|\rtil\omtil|^2-1+i\eps} \,d\omtil\rtil^2\,d\rtil\,d\eta \nn \\
&\les& R^3 \int_{\R^3} \int_{|\eta|-R^{-1}}^{|\eta|+R^{-1}}
\int_{[S^2:|\omtil-\frac{\eta}{|\eta|}|<R^{-1}]} \frac{|\hat{f}(\eta)|}{||\eta|-1|}
\frac{|\overline{\hat{f}(\rtil\omtil)}|}{|\rtil-1|} \,d\omtil\,d\rtil\,d\eta \nn \\
&\les& R^2 \int_{\R^3} \frac{|\hat{f}(\eta)|}{||\eta|-1|} \int_{|\eta|-R^{-1}}^{|\eta|+R^{-1}}
\left(\int_{[S^2:|\omtil-\frac{\eta}{|\eta|}|<R^{-1}]}
|\hat{f}(\rtil\omtil)|^2 \,d\omtil\right)^{\frac12}\,\frac{d\rtil}{|\rtil-1|}\,d\eta \nn \\
&\les& R^2 \int_0^\infty \frac{dr}{|r-1|} \int_{r-R^{-1}}^{r+R^{-1}}  \frac{d\rtil}{|\rtil-1|}
\int_{S^2} |\hat{f}(r\omega)| \left(\int_{[S^2:|\omtil-\omega|<R^{-1}]}
|\hat{f}(\rtil\omtil)|^2 \,d\omtil\right)^{\frac12} \nn 
\eear
and therefore also 
\bear
\eqref{eq:odin} &\les& R^2 \int_0^\infty \frac{dr}{|r-1|} \int_{r-R^{-1}}^{r+R^{-1}}  \frac{d\rtil}{|\rtil-1|}
\left(\int_{S^2} |\hat{f}(r\omega)|^2\,d\omega\right)^{\frac12} \left(\int_{S^2}\int_{[S^2:|\omtil-\omega|<R^{-1}]} |\hat{f}(\rtil\omtil)|^2 \,d\omtil\,d\omega\right)^{\frac12} \nn \\
&\les& R \int_{\frac12}^2 \frac{dr}{|r-1|} \int_{r-R^{-1}}^{r+R^{-1}}  \frac{d\rtil}{|\rtil-1|}
|1-r|^{\gamma}|1-\rtil|^{\gamma} \|f\|_p^2 \nn \\
&\les& R^{1-2\gamma}\|f\|_p^2 = R^{\frac{4}{p'}} \|f\|_p^2, \nn
\eear
where the last two lines use \eqref{eq:hold}. The lemma now follows by summing over
dyadic~$R$, at least provided~\eqref{eq:supp_res} holds. Finally, if
\[ \supp(\hat{f})\subset \{\xi\in\R^3\::\:|\xi|\le\half \text{\ or\ }|\xi|\ge2\}, \]
then one notes that
\[ \sup_{\eps\ne0}\|R_0(1\pm i\eps)f\|_2 \les \|(1-\triangle)^{-1}f\|_2 \les \|f\|_p \]
by the Sobolev imbedding theorem provided $1\le p\le 2$ and we are done.
\end{proof}

In Section~\ref{sec:more} we discuss further bounds on the free resolvent
which are motivated by the previous proposition.

\section{The perturbed resolvent}\label{sec:perturb}

The goal of this section is to prove theorem~\ref{thm:agmon}.
As in~\cite{Ag}, the proof of Theorem~\ref{thm:agmon} is based on the resolvent
identity. This requires inverting the operator $I+R_0(\lambda^2\pm i0)V$
on~$L^4(\R^3)$.  First, we check that this is a compact perturbation of the
identity.

\begin{lemma}
\label{lem:compact}
Let $V\in L^p(\R^3)$, $\frac32 \le p \le 2$. Then for any nonzero
$\lambda \in \Half$, the map $A(\lambda):=R_0(\lambda^2)V$
is a compact operator on $L^4(\R^3)$.
\end{lemma}
\begin{proof} Firstly, note that in view of Lemma~\ref{lem:free2} and because
of $V\in L^p$, $A(\lambda)$ is bounded $L^4\to L^4$.
Secondly, observe that we may assume that $V\in L^\infty$ with compact support.
Indeed, replace $V$
with $V_n=V\chi_{[|V|<n]}\chi_{[|x|<n]}$.
Then $\|V-V_n\|_p\to0$ as $n\to\infty$
implies that $\|A(\lambda)-A_n\|_{4\to 4}\to 0$ as $n\to\infty$.
If we can show that $A_n:=R_0(\lambda^2)V_n$ are
compact as operators $L^4\to L^4$ for each $n$, it therefore follows that
$A(\lambda)$ is also compact.
So assume that $V$ is bounded, and supported in the ball $\{|x| < R\}$.
Fix $\lambda$ and write $A=A(\lambda)$.
We first claim that $A:L^4\to W^{2,4}$. This follows from
\beeq
\label{eq:Aiden}
 (-\Laplace+1)A=(-\Laplace-\lambda^2)A+(\lambda^2+1)A = V+ (1+\lambda^2)A
\eneq
is bounded from $L^4$ to $L^4$.
Meanwhile, for $|x| > 2R$ there is the uniform pointwise bound
$$|Af(x)| \les \|Vf\|_1 |x|^{-1} \les
  R^{\frac94} \|V\|_\infty \|f\|_4 |x|^{-1}$$
Given $\eps > 0$, we may choose $R_0 \sim R^9\|V\|_\infty^4 \eps^{-4}$ so that
$\|\chi_{[|x|>R_0]} Af\|_4 < \eps$ for all $\|f\|_4 \le 1$.

Let $\{f_j\}_{j=1}^\infty\subset L^4(\R^3)$ satisfy
$f_j\rightharpoonup 0$ in $L^4$.  Since $\sup_j\|Af_j\|_{W^{2,4}(\R^3)}
< \infty$, Rellich's compactness theorem produces a subsequence $f_{j_k}$
so that $Af_{j_k} \to 0$ in $L^4(|x| < R_0)$.
Thus
\[ \limsup_{k\to\infty} \|Af_{j_k}\|_4 \le (1+C_\lambda)\eps.\]
Sending $\eps\to0$ and passing to the diagonal subsequence finishes the proof.
\end{proof}

The following lemma establishes invertibility everywhere except on the
imaginary axis.

\begin{lemma}
\label{lem:ev}
Let $V\in L^p(\R^3)\cap L^{\frac32}(\R^3)$, $\frac32 < p < 2$ and assume
that $V$ is real-valued.  Then for any nonzero $\lambda\in\Half$, the inverse
$(I+R_0(\lambda^2)V)^{-1}:L^4(\R^3)\to L^4(\R^3)$ exists.
\end{lemma}
\begin{proof}
By the previous lemma it suffices to show that
\[ f\in L^4(\R^3), \quad f+ R_0(\lambda^2)Vf=0 \Longrightarrow f=0.\]
Let $f$ be as on the left-hand side and set $g=Vf$. Then $g\in L^r$, where
$r= \frac{4p}{4+p} < \frac43$.
By Lemma~\ref{lem:free2}, $f = -R_0(\lambda^2)g$ therefore belongs to
$L^q \cap L^4$, where $\frac1q - \frac14 = \frac{3-2p}{3p} > 0$.

This bootstrapping procedure can be repeated until it is shown that $f \in
L^{r'} \cap L^4$.  In fact, one can continue to the point where
$f \in L^\infty$, since $R_0(\lambda^2): L^{\frac32 - \eps}\cap
L^{\frac32 + \eps} \mapsto L^\infty$ is a bounded operator.  What is important
here is that $f$ and $g$ exist in spaces dual to each other.

Since $V$ is real-valued, the duality pairing
\[ \langle f,g \rangle = \langle f, Vf \rangle = - \langle R_0(\lambda^2)g, g \rangle \]
shows that $\langle R_0(\lambda^2\pm i0)g, g \rangle$ is real-valued.
If $\lambda^2 \not\in \R$, then the condition
$$\Im \la R_0(\lambda^2)g, g\ra = \int_{\R^3} \frac{\Im(\lambda^2)}{\big(
 |\xi|^2 - \Re(\lambda^2)\big)^2 + \Im(\lambda^2)^2} |\hat{g}(\xi)|^2\,d\xi
 = 0 $$
requires that $\hat{g} = 0$ almost everywhere.

On the boundary $\lambda \in \R$, by the Stein-Tomas theorem
\[ \Im\langle R_0((\lambda + i0)^2)g, g \rangle =
\lim_{\eps\to 0}\Im \langle R_0((\lambda + i\eps)^2)g, g \rangle
= c\lambda \int_{S^2} |\hat{g}(\lambda\omega)|^2\, \sigma(d\omega)
\]
with some constant $c\ne0$. Hence, $\hat{g}=0$ on $|\lambda|S^2$ in the
$L^2$ sense.  Since $g\in L^r(\R^3)$, one concludes from
Proposition~\ref{prop:ess} above that
$(1+|x|)^{\delta-\half} R_0(\lambda^2\pm i0)g\in L^2(\R^3)$ for some
$\delta>0$.
Hence also $(1+|x|)^{\delta-\half}f \in L^2(\R^3)$ for some $\delta>0$.
Since $(-\Laplace+V-\lambda^2)f=0$ in the distributional sense, and one checks
easily from~\eqref{eq:Aiden} (remembering that $f \in L^\infty \cap L^4$)
that also $f\in W^{2,p}_{\rm loc}(\R^3) \subset W^{1,2}_{\rm loc}(\R^3)$,
Theorem~\ref{thm:IJ}  implies that $f=0$, as claimed.
\end{proof}

The following two lemmas show that the inverses in the previous lemma
have uniformly bounded norms.

\begin{lemma} \label{lem:continuous}
Let $V\in L^p(\R^3), \frac32 \le p \le 2$.  The map $\lambda \mapsto
R_0(\lambda^2)V$ is continuous from the domain $\Half \setminus \{0\} \subset
\Compl$ to the space of bounded operators on $L^4(\R^3)$.
\end{lemma}
\begin{proof}
First suppose $V$ is bounded and has compact support in the ball $\{|x| < R\}$.
The convolution kernel associated to $R_0(\lambda^2) - R_0(\zeta^2)$ has the
bounds
\[ |K(x)| \les \left\{
\begin{aligned} |\lambda-\zeta|, \quad &{\text if}\ |x|<|\lambda-\zeta|^{-1} \\
                |x|^{-1}, \quad &{\text if}\ |x| \ge |\lambda-\zeta|^{-1}
\end{aligned}\right.
\]
Then for any pair $\lambda,\zeta \in \Half$, $|\lambda-\zeta| \le \frac1{2R}$,
we have
\[|(R_0(\lambda^2)-R_0(\zeta^2))Vf(x)| \les \left\{
\begin{aligned} |\lambda-\zeta|\norm[Vf][1], \quad &{\text if}\
   |x| < |\lambda-\zeta|^{-1} \\
   |x|^{-1}\norm[Vf][1], \quad &{\text if}\ |x| \ge |\lambda-\zeta|^{-1}
\end{aligned} \right.
\]
Thus $\norm[(R_0(\lambda^2)-R_0(\zeta)^2)Vf][4] \les |\lambda-\zeta|^{1/4}
R^{9/4}\norm[V][\infty]\norm[f][4]$.

Approximate $V$ by compactly supported $\tilde{V}\in L^\infty$ so that
$\norm[V-\tilde{V}][p] < \eps$.  By the above calculation,
Lemma~\ref{lem:free2}, and the simple identity
$$(R_0(\lambda^2)-R_0(\zeta)^2)V = R_0(\lambda^2)(V-\tilde{V})
 + (R_0(\lambda^2) - R_0(\zeta^2))\tilde{V} - R_0(\zeta^2)(V-\tilde{V}),$$
we see that $\limsup_{\zeta\to\lambda}
\norm[(R_0(\lambda^2)-R_0(\zeta)^2)V][4\to 4] \les |\lambda|^{(3-2p)/p}\eps$.
\end{proof}

\begin{lemma}
\label{lem:unif}
Let $V$ be as in the previous lemma and suppose $\lambda_0>0$. Then
\beeq
\label{eq:unif}
\sup_{|\Re(\lambda)|\ge\lambda_0}
\Big\| (I+R_0(\lambda^2)V)^{-1} \Big\|_{4\to 4} < \infty.
\eneq
\end{lemma}
\begin{proof}
In view of Lemma~\ref{lem:free2}, there is some finite $\lambda_1 \in \R$
so that $\norm[R_0(\lambda^2)V][4\to4] < \frac12$ provided
$|\lambda| > \lambda_1$.  It therefore suffices to prove \eqref{eq:unif} on the
compact set $\{\lambda\in\Compl:\;\lambda_0 \le |\lambda| \le \lambda_1,
|\Re(\lambda)| \ge \lambda_0\}$.  The previous two lemmas, however, show that
$(I + R_0(\lambda^2)V)^{-1}$ is a continuous function of $\lambda$ on this
set, hence it is uniformly bounded from above.
\end{proof}

It is now a simple matter to prove Theorem~\ref{thm:agmon}.

\begin{proof}[Proof of Theorem~\ref{thm:agmon}.]
By the resolvent identity, for any $\eps\ne0$,
\[ R_V(\lambda^2+i\eps)=R_0(\lambda^2+i\eps)-R_0(\lambda^2+i\eps)VR_V(\lambda^2+i\eps).\]
By Lemma~\ref{lem:unif} one therefore has
\[ R_V(\lambda^2+i\eps)= (I+R_0(\lambda^2+i\eps)V)^{-1}R_0(\lambda^2+i\eps)\]
and the right-hand side is uniformly bounded for $\lambda\ge\lambda_0\ge0$ as well
as $0<\eps\le1$ in the $L^4$ operator norm.
In fact, the last factor contributes a decaying factor of $\lambda^{-\half}$
as $L^4$ operator norm in view of Lemma~\ref{lem:free}.
\end{proof}

\begin{proof}[Proof of Proposition~\ref{prop:extension}.]
There is only one point in the argument where the condition
$V\in L^{\frac32}(\R^3)$ is used, namely the step in Lemma~\ref{lem:ev} where
we wish to make use of Theorem~\ref{thm:IJ}.  It otherwise suffices to assume
that $V \in L^p(\R^3)$, $\frac32 < p < 2$.

For the second claim, one observes the following consequence of
Lemma~\ref{lem:free2}: If $V \in L^p$, $\frac32 < p \le 2$, and $r > 4$, then
$R_0(\lambda^2\pm 10)V: L^4 \cap L^r \mapsto L^4 \cap L^s$, where
$\frac1s = \max(\frac1r + \frac1p - \frac23, 0)$.  The same is true
for any $V\in L^q, p \le q \le 2$.  This allows the
bootstrapping procedure on $f$ to continue normally, and furthermore
$g = Vf$ is still an element of $L^{\frac43-\eps}$, as desired.
Therefore, the only matter of concern is whether the conclusion of
Theorem~\ref{thm:IJ} will hold for such a potential $V$.
\end{proof}

\section{Further estimates on the free resolvent}\label{sec:more}

Returning to Proposition~\ref{prop:ess}, we note that a
sharper estimate can be made at the endpoint $p=1$.

\begin{proposition}\label{prop:p=1}
Let $f$ be a function in $L^1(\R^3)$ such that $\hat{f} = 0$ on the unit sphere
$S^2$.  Then
\begin{equation}
\sup_{\eps> 0} \big\|R_0(1\pm i\eps)f\big\|_2 \le \frac1{\sqrt{8\pi}}
  \norm[f][1]
\end{equation}
\end{proposition}
\begin{proof}
Define the trace function
\begin{equation} \label {eq:G1}
G(\lambda) = \lambda^{-2}\big\|\hat{f}|_{\lambda S^2}\big\|_2^2
  = 4\pi\iint_{\R^6} f(x) \frac{\sin(\lambda|x-y|)}{\lambda|x-y|} \bar{f}(y)
\,dx\,dy
\end{equation}
By inspection,
\begin{equation} \label{eq:G2}
G(\lambda) = 2\pi \iiint_{\R\times\R^6}
 \frac{f(x)\bar{f}(y)}{|x-y|} \chi_{|x-y|}(\tau) e^{i\lambda\tau}\,
  d\tau dx dy
\end{equation}
where $\chi_{|x-y|}$ denotes the characteristic function of the interval
 $\{|\tau| \le |x-y|\}$.  The integrand on the right-hand side is in
$L^1(\R^7)$, so Fubini's Theorem implies that $G$ is the inverse Fourier
transform of an $L^1$ function.

Using the Plancherel identity (in 3 dimensions), and noting that $G$ is an even
function,
\begin{equation} \label{eq:GM}
\norm[R_0(1\pm i\eps)f][2]^2 =  \frac1{2(2\pi)^3}\int_{-\infty}^\infty
  G(\lambda)\frac{\lambda^2}{|\lambda^2- (1+i\eps)|^2}
\end{equation}
For any $\eps > 0$, the multiplier $M_{\eps}(\lambda) =
\frac{\lambda^2}{|\lambda^2 - (1+i\eps)|^2}$ is integrable, hence it has
Fourier transform $\hat{M}_\eps \in L^\infty(d\tau)$.  By Parseval's
formula, this time in one dimension,
\begin{equation} \label{eq:Mepsilon}
\norm[R_0(1\pm i\eps)f][2]^2 = \frac1{2(2\pi)^4}\int_\R\hat{G}(\tau)
 \hat{M}_\eps(-\tau) \,d\tau
\end{equation}
An explicit formula for $\hat{M}_\eps(\tau)$ can be obtained via
residue integrals:
\begin{equation}
\hat{M}_\eps(\tau) =
\frac\pi{2\eps}\big(\sqrt{1+i\eps}\,e^{i|\tau|\sqrt{1+i\eps}}
                       +\sqrt{1-i\eps}\,e^{-i|\tau|\sqrt{1-i\eps}}\big)
\end{equation}
This, along with \eqref{eq:G2}, can be immediately substituted back into
equation \eqref{eq:Mepsilon}.
\begin{equation*}
\begin{aligned}
\norm[R_0(1\pm i\eps)f][2]^2 &=
\frac{1}{8\pi\eps}\iint_{\R^6}\int_0^{|x-y|} \frac{f(x)\bar{f}(y)}{|x-y|}
\big(\sqrt{1+i\eps}\,e^{i\tau\sqrt{1+i\eps}}
      +\sqrt{1-i\eps}\,e^{-i\tau\sqrt{1-i\eps}}\big)\,d\tau\,dx\,dy \\
&= \frac{1}{8\pi i\eps}\iint_{\R^6}\frac{f(x)\bar{f}(y)}{|x-y|}
   \big(e^{i|x-y|\sqrt{1+i\eps}} - e^{-i|x-y|\sqrt{1-i\eps}}\big)dx\,dy
\end{aligned}
\end{equation*}
Boundedness of $\hat{M}_\eps$ enables us to continue applying Fubini's
theorem to the multiple integral.  We have also simplified the expression by
noting that $\hat{M}_\eps$ is an even function.  Recall definition
\eqref{eq:G1} and subtract $\frac1{16\pi^2\eps}G(1)$ from both sides of
the equation.
\begin{equation}
\begin{aligned}
\|R_0(1&\pm i\eps)f\|_2^2 - \frac1{16\pi^2\eps}G(1) \\
 &= \frac{1}{8\pi i\eps}\iint_{\R^6}\frac{f(x)\bar{f}(y)}{|x-y|}
   \Big(\big(e^{i|x-y|\sqrt{1+i\eps}}-e^{i|x-y|}\big) -
        \big(e^{-i|x-y|\sqrt{1-i\eps}}-e^{-i|x-y|}\big)\Big) dx\,dy \\
 &= \frac1{8\pi i\eps}\iint_{\R^6}\frac{f(x)\bar{f}(y)}{|x-y|}K(|x-y|)
\,dx\,dy
\end{aligned}
\end{equation}
where $|K(|x-y|)| \le \eps|x-y|$.  This leads to the conclusion
$$\big| \norm[R_0(1\pm i\eps)f][2]^2 - \frac1{16\pi^2\eps}G(1)\big|
  \le \frac1{8\pi} \norm[f][1]^2$$
If $f$ satisfies the hypothesis $\hat{f}|_{S^2} = 0$, then $G(1) = 0$.
\end{proof}

\begin{cor}
Let $f$ be a function in $L^1(\R^3)$ such that $\hat{f} = 0$ on the unit sphere
$S^2$.  Then
\begin{equation}
 \norm[R_0(1\pm i0)f][2] \le \frac1{\sqrt{8\pi}} \norm[f][1]
\end{equation}
\end{cor}
\begin{proof}
This follows immediately from \eqref{eq:GM} and monotone convergence.
\end{proof}

The condition $\hat{f}=0$ is crucial in Proposition~\ref{prop:ess}.
Indeed, recall that for $f\in L^p(\R^3)$ real-valued with $1\le p\le\frac43$ one has
\[ \Im R_0(1+i0)f = c\,(\widehat{\sigma_{S^2}}\ast f) \]
for some constant $c$. This follows by writing $R_0(1+i\eps)$ as a sum of its real
and imaginary parts, as well as from the fact that the operation of restriction
$f\mapsto \hat{f}(r\cdot)$ is continuous in $r>0$ as a map $L^p(\R^3)\to L^2(S^2)$.
However, it is clear that for any $\delta>0$
\begin{equation}
\label{eq:inf}
 \|(1+|x|)^{\delta-\half} [\widehat{\sigma_{S^2}}\ast f] \|_{2} = \infty
\end{equation}
even for smooth bump-functions $f$ since the function inside the norm decays like
$(1+|x|)^{\delta-\frac32}$ which just fails to be $L^2(\R^3)$.
The following simple lemma shows, on the other hand, that $\delta<0$ does lead to a finite
norm in~\eqref{eq:inf}.

\begin{lemma}
\label{lem:localST}
For any $R\ge 1$ one has
\[
\Big\|\chi_{[|x|<R]}  [\widehat{\sigma_{S^2}}\ast f] \Big\|_{2} \les \sqrt{R}\; \|f\|_{\frac43}
\]
for all $f\in L^{\frac43}(\R^3)$.
\end{lemma}
\begin{proof}
Let $\phi$ be a smooth cut-off function with $\hat{\phi}$ compactly supported.
Then by Plancherel, and Cauchy-Schwartz,
\bear
&& \Big\|\chi\Big(\frac{\cdot}{R}\Big) [\widehat{\sigma_{S^2}}\ast f] \Big\|^2_{2}
= R^6 \int_{\R^3} \left| \int_{S^2} \hat{\chi}(R(\xi-\eta)) \hat{f}(\eta)\,\sigma_{S^2}(d\eta)
\right|^2 \, d\xi \nn \\
&& \les R^6 \int_{\R^3} \int_{S^2} |\hat{\chi}(R(\xi-\eta'))| \,d\eta' \int_{S^2}
|\hat{\chi}(R(\xi-\eta))| |\hat{f}(\eta)|^2\,\sigma_{S^2}(d\eta)\,d\xi \nn \\
&&\les R\; \|\hat{f}\|^2_{L^2(S^2)} \les R\; \|f\|^2_{\frac43}, \nn
\eear
as claimed.
\end{proof}

The previous lemma suggests that one should also have the bound
\begin{equation}
\label{eq:thorn}
\sup_{\eps>0}\Big\|\chi_{[|x|<R]} R_0(1\pm i\eps) f \Big\|_{2} \les \sqrt{R}\; \|f\|_{\frac43}.
\end{equation}
While this bounds remains open\footnote{Note added in proof: This conjecture
is solved in a forthcoming paper by Ionescu and the second author~\cite{IonSch}, which also contains 
other improvements on the results obtained here.}, it is easy to show that
\begin{equation}
\label{eq:34}
\sup_{\eps>0}\Big\|\chi_{[|x|<R]} R_0(1\pm i\eps) f \Big\|_{2} \les R^{\frac34}\; \|f\|_{\frac43}.
\end{equation}
Indeed, denoting the operator on the left-hand side by $T$ for a fixed $\eps>0$, observe that
by Lemma~\ref{lem:free}
\[ T^*T = R_0(1- i\eps)\chi_{[|x|<R]} R_0(1+ i\eps) \]
satisfies $\|T^*T\|_{\frac43\to 4} \les \| \chi_{[|x|<R]} \|_2 \les R^{\frac32}$, which is the
same as~\eqref{eq:34}. One would of course expect that the Knapp example determines the power
of~$R$ in~\eqref{eq:thorn}. Note that in the case of a Knapp example of dimensions
$\delta\times\sqrt{\delta}\times\sqrt{\delta}$ where $\delta=R^{-1}$, one does not encounter the
$L^2$ norm of $\| \chi_{[|x|<R]} \|_2$ in the previous $T^*T$ argument, but rather the $L^2$ norm
of a $R\times \sqrt{R}\times\sqrt{R}$-tube, which gives the conjectured~$R$.
Conversely, in what follows we show how to approach~\eqref{eq:thorn} by a decomposition
into Knapp examples. This leads to an improvement over the simple bound~\eqref{eq:34} by~$\frac18$.
Such arguments originate in the analysis of Bochner Riesz
multipliers as well as the restriction theory of the Fourier transform.
We will use a square function bound from Bourgain~\cite{B}. For the convenience
of the reader, we reproduce the details.

The following lemma is a discrete version of the Stein-Tomas
theorem, see Lemma~6.2 in~\cite{B}. It can also be proved by
expanding the $L^4$-norm on the left-hand side explicitly and then
using the usual geometric arguments based on counting overlap.
However, the following approach does not depend on any special
arithmetic properties of the Stein-Tomas exponent and therefore
generalizes to other dimensions as well.

\begin{lemma}
\label{lem:jean}
Let $R\ge1$ and let $\{\xial\}_\alpha\subset S^2$ be a collection of $R^{-\half}$-separated points and
$\{a_\alpha\}_\alpha$ a sequence of arbitrary numbers. Then
\begin{equation}
\label{eq:disc}
\Big\|\sum_\alpha e^{i x\cdot \xial}\, a_\alpha \Big\|_{L^4(Q)} \les R^{\half} \Big(\sum_\alpha |a_\alpha|^2 \Big)^{\half}
\end{equation}
for any cube $Q$ of side-length $R^{\half}$.
\end{lemma}
\begin{proof}
Fix $R\ge1$, some cube $Q$ and a smooth cut-off function $\chi_Q$ adapted
to~$Q$. We assume that $\supp(\widehat{\chi_Q})$ is contained inside a $CR^{-\half}$-cube
with some big constant $C$.
Then
\bear
\Big\| \sum_\alpha e^{ix\cdot\xial} a_\alpha \Big\|_{L^4(Q)} &\les&
\Big\| \Big[\sum_\alpha a_\alpha \widehat{\chi_Q}(\cdot-\xial)\Big]^{\vee}\Big\|_4 \nn \\
&\les& \int_{1-R^{-\half}}^{1+R^{-\half}} \Big\| \int_{S^2} e^{ir\omega\cdot x}\;\sum_\alpha
a_\alpha \widehat{\chi_Q}(r\omega-\xial) \; \sigma(d\omega) \Big\|_{L^4_x} \, dr \nn \\
&\les& R^{-\half} R^{\frac32} \Big(\sum_\alpha R^{-1}\,|a_\alpha|^2 \Big)^{\half}
= R^{\half}\,\Big(\sum_\alpha |a_\alpha|^2 \Big)^{\half}, \nn
\eear
as claimed.
\end{proof}

In what follows, we associate with every $R\ge1$ a decomposition $\{\Sal\}_{\alpha}$ of the shell
$\calS_R:=\{\xi\in \R^3\::\: ||\xi|-1|<R^{-1}\}$
into Knapp caps  $\Sal$ of size about $R^{-\half}\times R^{-\half}\times R^{-1}$.
We furthermore assume bounded overlap, i.e., that
\[ \sup_{R\ge1}\Big\|\sum_\alpha \chi_{\Sal}\Big\|_\infty  < \infty.\]

\begin{lemma}
\label{lem:square}
Let $g\in L^4(\R^3)\cap L^2(\R^3)$ be so that $\hat{g}\subset \calS_R$.
Write $g=\sum_\alpha \gal$
where each $\widehat{\gal}$ is supported in~$\Sal$. Then
\begin{equation}
\label{eq:g2}
\|g\|_4 \les R^{\frac18} \Big\|\Big(\sum_\alpha \big|\gal\big|^2\Big)^{\frac12}\Big\|_4.
\end{equation}
Dually, one has for any $f\in L^{\frac43}(\R^3)$ that
\begin{equation}
\label{eq:f2}
\Big\|\Big(\sum_\alpha \big|\fal\big|^2\Big)^{\frac12}\Big\|_{\frac43} \les R^{\frac18}\; \|f\|_{\frac43}
\end{equation}
where $\widehat{\sum_\alpha \fal}=\hat{f}\chi_{\calS_R}$ and each $\fal$ is Fourier supported
in $\Sal$.
\end{lemma}
\begin{proof}
Firstly, fix $R\ge1$, a cube $Q$ of size $\sqrt{R}$, the cap
decomposition $\{\Sal\}_\alpha$, as well as some $\xial\in\Sal$
for each $\alpha$. Secondly, fix some smooth cut-off $\chi_0$ so
that $\chi_0(\sqrt{R}\,y)$ is adapted to the cube
$Q_0:=[-c_0\sqrt{R},c_0\sqrt{R}]^3$ where $c_0>0$ is small. We
require that $\chi_0(\sqrt{R}\,(\xi-\xial))=1$ for all
$\xi\in\Sal$ and each $\alpha$. Set \[
\mu(\xi):=\chi_0(\sqrt{R}\,\xi)\left(\slashint_{Q_0} e^{i
y\cdot\xi} \, dy\right)^{-1}\] and
$\mu_\alpha(\xi)=\mu(\xi-\xial)$ for every $\alpha$. Note that
this is again a smooth cut-off function adapted to a ball of size
$R^{-\half}$ (together with the natural derivative bounds with
constants uniform in~$R$). In particular, $\sup_\alpha
\|\widehat{\mu_\alpha}\|_1 = \|\hat{\mu}\|_1 < \infty$ uniformly
in~$R$. Then \bear \big\|g\big\|_{L^4(Q)}^4 &=&
\Big\|\sum_\alpha \int e^{ix\cdot\xi}\; \widehat{g_\alpha}(\xi)\, d\xi \Big\|_{L^4(Q)}^4 \nn \\
&\le& \int_{Q}\slashint_{Q_0}
\left|\sum_\alpha e^{-iy\cdot\xial} \int e^{i(x+y)\cdot\xi}\; \widehat{g_\alpha}(\xi)\mu_\alpha(\xi)\, d\xi \right|^4\,dy\,dx \label{eq:jens} \\
&\les& \slashint_{2Q} \int_{Q_0}  \left|\sum_\alpha e^{-iy\cdot\xial}\; \int e^{ix\cdot\xi} \widehat{g_\alpha}(\xi)\mu_\alpha(\xi)\, d\xi \right|^4\,dy\,dx \label{eq:change} \\
&\les& \slashint_{2Q} R^2  \left(\sum_\alpha \left|\int e^{ix\cdot\xi}\; \widehat{g_\alpha}(\xi)\mu_\alpha(\xi)\, d\xi \right|^2 \right)^2 \,dx \label{eq:bou} \\
&\les& R^{\half} \Big\|\Big(\sum_\alpha |h_\alpha|^2\Big)^{\half}\Big\|^4_{L^4(2Q)}, \label{eq:hal}
\eear
where $h_\alpha=g_\alpha\ast k_\alpha$ and $\widehat{k_\alpha}=\mu_\alpha$.
Here \eqref{eq:jens} follows from Jensen's inequality as well as the definition of $\mu_\alpha$,
\eqref{eq:change} follows by changing variables (and $2Q$ is the cube with the same center
but twice the size), whereas~\eqref{eq:bou} is a consequence of
Lemma~\ref{lem:jean}. Summing~\eqref{eq:hal} over
a partition $\{Q\}$ of $\R^3$ consisting of congruent cubes one obtains that
\begin{align*}
 \big\|g\big\|_{4} & \les R^{\frac18} \Big\|\Big(\sum_\alpha |h_\alpha|^2\Big)^{\half}\Big\|_4
\les R^{\frac18} \Big\|\Big(\sum_\alpha |g_\alpha|^2\ast\sup_{\alpha}|k_\alpha|\Big)^{\half}\Big\|_4
\les R^{\frac18} \Big\|\sum_\alpha |g_\alpha|^2 \ast |\check{\mu}| \Big\|_2^{\half} \\
& \les R^{\frac18} \Big\| \Big(\sum_\alpha |g_\alpha|^2 \Big)^{\half} \Big\|_4
\end{align*}
since $ \|\check{\mu}\|_1 \les 1$.
This proves~\eqref{eq:g2} and~\eqref{eq:f2} follows by duality. Indeed,
\bear
\|\{f_\alpha\}\|_{L^{\frac43}(\ell^2)} &=& \sup_{\|\{g_\alpha\}\|_{L^4(\ell^2)}\le 1}
 \left|\int \sum_\alpha f_\alpha\; \overline{g_\alpha} \right| \les
\sup_{\|g\|_4\le CR^{\frac18},\;{\rm supp}(\hat{g})\subset \calS_R} |\la f,g\ra | \label{eq:geins}\\
&\les& \sup_{\|g\|_4\le CR^{\frac18}} |\la f,g\ra | \les R^{\frac18} \|f\|_{\frac43}.\label{eq:gzwei}
\eear
In line~\eqref{eq:geins}, the supremum is taken over all $g$ as in~\eqref{eq:g2}, whereas
in~\eqref{eq:gzwei}, we drop the condition $\supp(\hat{g})\subset \calS_R$.
\end{proof}

The powers of $R$ in~\eqref{eq:disc}, as well as
in~\eqref{eq:g2} and~\eqref{eq:f2} are optimal. In the case of~\eqref{eq:disc}, this can be
seen by taking all $a_\alpha=1$ and similarly for the square function.

\noindent
Recall that Agmon's
limiting absorption
principle states that
\begin{equation}
\label{eq:limap}
\sup_{\eps>0}\big\|(1+|x|)^{-\half-}R_0(1\pm i\eps) f\big\|_2 \les \big\|(1+|x|)^{\half+}f\big\|_2.
\end{equation}
This in particular implies that
\begin{equation}
\label{eq:weaklimap}
\sup_{\eps>0} \big\|\chi_{[|x|<R]} R_0(1\pm i\eps)\chi_{[|x|<R]} f \Big\|_{2} \les R^{1+}\;
\|f\|_{2}.
\end{equation}
The following lemma is an improved version of~\eqref{eq:weaklimap},
see the decay in~$|v|$. The proof is self-contained. In particular, it does
not rely on~\eqref{eq:limap}.

\begin{lemma}
\label{lem:limap2}
For any $v\in\R^3$ one has
\begin{equation}
\label{eq:L2limap}
\sup_{\eps>0}\Big\|\chi_{[|x|<R]} R_0(1\pm i\eps)\chi_{[|x-v|<R]} f \Big\|_{2} \les R(1+|v|/R)^{-1}\;
\|f\|_{2}.
\end{equation}
for all $R\ge1$. The constants are uniform in $R$ and $v$.
\end{lemma}
\begin{proof}
This is basically a simple consequence of  H\"ormander's variable
coefficient Plancherel theorem, see~\cite{hor}, and especially Wolff's notes~\cite{wol}, page~55.
We start with the case $|v|\gg R$. Consider the operator
\[ T_{R,v}f(x) = \int_{\R^3}\psi(x)\frac{e^{iR|x-y|}}{|x-y|}\psi(y-v/R)f(y)\,dy \]
where $\psi$ is a  smooth bump function at zero. This operator is of the form
\[ T_{R,v}f(x)= \int_{\R^3} e^{iR\Phi(x,y)}a_v(x,y) f(y)\,dy \]
where $\Phi$ is smooth on the support of $a_v(x,y)$, and $\rank[ \partial^2_{xy} \Phi(x,y)]=2$.
Moreover, $a_v(x,y)$ is smooth, the size of its support is uniformly bounded in $v$,  and
$\|\partial^{\beta} a_v\|_\infty\les (1+|v|/R)^{-1}$. Hence, H\"ormander's variable coefficient Plancherel
theorem  implies that
\[ \|T_{R,v}\|_{2\to2} \les R^{-1}(1+|v|/R)^{-1}.\]
One now checks that~\eqref{eq:L2limap} follows from this by means of a change of variables.

\noindent If $|v|\les R$, one can argue similarly, but needs to introduce a Whitney
decomposition away from the singularity $x=y$. Firstly, note that with $\psi$ as above,
\[
\left\|\int_{\R^3}\psi(x)\frac{e^{iR|x-y|}}{|x-y|}\chi_{[|x-y|\le R^{-\half}]}\psi(y)f(y)\,dy \right\|_2
\les \int_{[|x|< R^{-\half}]} \frac{dx}{|x|} \; \|f\|_2 \les R^{-1}\,\|f\|_2.
\]
Next, let $\rho\in C^\infty_0(\R)$ be a smooth function so that $\rho(t)=1$ if
$1<t<2$ and $\rho(t)=0$ if $t>2$ or $t<\half$. Then we claim that
\beeq
\label{eq:whitney}
\Big\|\int_{\R^3}\psi(x)\frac{e^{iR|x-y|}}{|x-y|}\rho(2^{-j}|x-y|)\psi(y)f(y)\,dy \Big\|_2
\les 2^j\,R^{-1}\,\|f\|_2
\eneq
for all $R^{-\half}<2^j\le 1$. To this end introduce a further decomposition
$1=\sum_\ell \omega\big(2^{-j}(x-x_\ell)\big)$ where the sum runs over a
lattice of points $\{x_\ell\}_\ell$ in~$\R^3$ that are $2^j$-spaced, and $\omega$
is some smooth cut-off which is adapted to the unit cube.
 Exploiting orthogonality, \eqref{eq:whitney} follows from
the following estimate
\[
\Big\|\int_{\R^3}\omega\big(2^{-j}(x-x_\ell)\big)
\frac{e^{iR|x-y|}}{|x-y|}\rho(2^{-j}|x-y|)\omega\big(2^{-j}(y-x_k)\big)\;f(y)\,dy \Big\|_2
\les 2^j\,R^{-1}\,\|f\|_2.
\]
This, however, is again reduced the H\"ormander's bound by means of an obvious rescaling.
\end{proof}

We are now ready to formulate our estimate which lies between the conjecture~\eqref{eq:thorn}
and the simple bound~\eqref{eq:34}.

\begin{prop}
\label{prop:58}
Let $f\in L^{\frac43}(\R^3)$.
Then there is the bound
\begin{equation}
\label{eq:58}
\sup_{\eps>0}\big\|\chi_{[|x|<R]} R_0(1\pm i\eps) f \big\|_{2} \les R^{\frac58}\;
\|f\|_{\frac43}.
\end{equation}
for all $R\ge1$.
\end{prop}
\begin{proof}
If $\supp(\hat{f})\subset \big\{\xi\in\R^3\::\:||\xi|-1|>\half\big\}$, then
one has
\[ \sup_{\eps>0}\big\|\chi_{[|x|<R]} R_0(1\pm i\eps) f \big\|_{2} \les
\|(-\Laplace+1)^{-1}f\|_2 \les \|f\|_{\frac43}.\]
Hence we can assume that $\supp(\hat{f})\subset \big\{\xi\in\R^3\::\:||\xi|-1|\le \half\big\}$
and we write $f=\sum_k f_k$, where
\[ \supp(\hat{f_k})\subset \big\{\xi\in\R^3\::\:||\xi|-1|\asymp 2^k R^{-1}\big\} \]
for $k\ge 1$ and
\[ \supp(\hat{f_0})\subset \big\{\xi\in\R^3\::\:||\xi|-1|\le R^{-1}\big\}. \]
Now $R_0(1+i\eps)f_0=R_0(1+i\eps)\sum_\alpha f_\alpha$
with $\sum_\alpha f_\alpha$ as in Lemma~\ref{lem:square}.
Consider any of the pieces $\fal$. Then one can write $\fal=\sum_j f_{\alpha,j}$ where
each $f_{\alpha,j}$ lives on a tube of dimensions $R\times\sqrt{R}\times\sqrt{R}$ and the tubes
corresponding to different $j$'s are disjoint. Moreover, the bound
\[ \|f_{\alpha,j}\|_\infty \les \|f_{\alpha,j}\|_p\, R^{-\frac{2}{p}} \]
holds for every $1\le p\le\infty$. Therefore,
\bear
\|\fal\|_2^2 &\les& \sum_j \|f_{\alpha,j}\|_2^2 \les \sum_j \|f_{\alpha,j}\|_\infty^2 R^2 \nn \\
&\les& \sum_j \|f_{\alpha,j}\|_{\frac43}^2 R^{-1} \les \left(\sum_j \|f_{\alpha,j}\|_{\frac43}^{\frac43}\right)^{\frac32}\; R^{-1} \nn \\
&\les& \|\fal\|_{\frac43}^2 \;R^{-1}. \label{eq:243}
\eear
We will apply this bound not to $\fal$, but to $\psi((\cdot-v)/R)\fal$ for a smooth cut-off $\psi$
and arbitrary $v$. This is justified since the Fourier support of $\psi((\cdot-v)/R)\fal$
still lies in the cap~$\Sal$.
In combination with~\eqref{eq:L2limap}, \eqref{eq:243} yields
\bear
&& \sup_{\eps>0} \big\|\chi_{[|x|<R]} R_0(1\pm i\eps) \fal \Big\|_{2}
\les \sum_{v\in R\Z^3}
\sup_{\eps>0} \big\|\chi_{[|x|<R]} R_0(1\pm i\eps)\chi_{[|x-v|<R]} \fal \Big\|_{2} \nn \\
&& \les R\; \sum_{v\in R\Z^3}  (1+|v|/R)^{-1}
\|\psi((\cdot-v)/R)\fal\|_{2} \nn \\
&& \les  R\; \sum_{v\in R\Z^3}  (1+|v|/R)^{-1} \, R^{-\half}\,
\|\psi((\cdot-v)/R)\fal\|_{\frac43} \nn \\
&& \les R^{\half} \Big(\sum_{j\in\Z^3}(1+|j|)^{-4}\Big)^{\frac14}\Big(\sum_{j\in\Z^3}\|\psi((\cdot-Rj)/R)\fal\|_{\frac43}^{\frac43}\Big)^{\frac34} \les R^{\half}\|f_\alpha\|_{\frac43}. \label{eq:knapp_agmon}
\eear
It remains to sum up~\eqref{eq:knapp_agmon} exploiting  orthogonality as provided by
Lemma~\ref{lem:square}:
\bear
\sup_{\eps>0}\Big\|\chi_{[|x|<R]} R_0(1\pm i\eps) f_0 \Big\|_{2}^2 &\le&
\sum_\alpha \sup_{\eps>0} \Big\|\psi\big(\frac{\cdot}{R}\big) R_0(1\pm i\eps) f_\alpha \Big\|_{2}^2 \nn \\
&\les& R \sum_\alpha \|  f_\alpha \|_{\frac43}^2 \les R^{\frac54} \|f\|_{\frac43}^2. \nn
\eear
The final bound uses \eqref{eq:f2}. It is easier to deal with $f_k$ where $k\ge1$.
Fix such a $k$ and let $f_k=\sum_\beta f_\beta$ be the corresponding decomposition
into Knapp caps. Then, by~\eqref{eq:243},
\[ \|f_\beta\|_2 \les \|f_\beta\|_{\frac43} (2^k\,R^{-1})^{\half}\]
and instead of~\eqref{eq:knapp_agmon} one now has
\[ \sup_{\eps>0} \big\|\chi_{[|x|<R]} R_0(1\pm i\eps) f_\beta \Big\|_{2}
\les (2^k R^{-1})^{-1} \|f_\beta\|_2 \les (2^{-k} R)^{\half} \|f_\beta\|_{\frac43}.\]
Invoking~\eqref{eq:f2} finally allows one to conclude that
\[ \sup_{\eps>0} \big\|\chi_{[|x|<R]} R_0(1\pm i\eps) f_k \Big\|_{2}
\les (2^{-k} R)^{\frac54} \|f\|_{\frac43}.\]
Summing over $k$ finishes the proof.
\end{proof}

{\bf Acknowledgement:} The second author was supported by the NSF grant
DMS-0300081 and a Sloan fellowship. The authors wish to thank the referee
for numerous suggestions.

\bibliographystyle{amsplain}

\medskip\noindent
\textsc{Division of Astronomy, Mathematics, and Physics, 253-37 Caltech, Pasadena, CA 91125, U.S.A.}\\
{\em email: }\textsf{\bf mikeg@caltech.edu, schlag@caltech.edu}

\end{document}